\documentclass{amsart}
\usepackage{amsmath, amsthm, amscd, amsfonts, amssymb, graphicx, color, enumerate}
\usepackage[bookmarksnumbered, colorlinks, plainpages]{hyperref}
\usepackage[all]{xy}

\newtheorem{theorem}{Theorem}[section]
\newtheorem{lemma}[theorem]{Lemma}
\newtheorem{proposition}[theorem]{Proposition}
\newtheorem{corollary}[theorem]{Corollary}
\theoremstyle{definition}
\newtheorem{definition}[theorem]{Definition}

\theoremstyle{remark}
\newtheorem{remark}[theorem]{Remark}
\numberwithin{equation}{section}

\newcommand{\V}[1]{\operatorname{V}\left(#1\right)}

\newcommand{\tor}[4]{\operatorname{Tor}^{#1}_{#2}\left(#3,#4\right)}
\newcommand{\h}[3]{\operatorname{H}^{#1}_{#2}\left(#3\right)}

\newcommand{\gh}[4]{\operatorname{H}^{#1}_{#2}(#3,#4)}

\newcommand{\ann}{\operatorname{Ann}}

\newcommand\supp {\operatorname{Supp}}

\newcommand\att{\operatorname{Att}}

\newcommand{\width}[2]{\operatorname{width}(#1,#2)}
\newcommand{\fwidth}[2]{\operatorname{f-width}(#1,#2)}

\newcommand\im{\operatorname{Im}}

\newcommand\Ker{\operatorname{Ker}}

\newcommand\cosupp {\operatorname{Cosupp}}
\newcommand\Max{\operatorname{Max}}

\newcommand\fa{\mathfrak a}

\newcommand\fp{\mathfrak p}
\newcommand\fq{\mathfrak q}

\newcommand\N{\mathbb N}
\swapnumbers

\begin{document}
\setcounter{page}{1}



\title[On the finiteness of local homology modules]{On the finiteness of local homology modules}
\author[Fathi and Hajikarimi]{A. Fathi$^{*}$ and A.  Hajikarimi}
\thanks{{\scriptsize
\hskip -0.4 true cm MSC(2010): Primary: 13D07; Secondary: 13E05, 13C05, 13C15, 13D45.
$*$Corresponding author }}

\begin{abstract}
Let $R$ be a commutative Noetherian ring and $\mathfrak{a}$ be an ideal of $R$. Suppose $M$ is a finitely generated $R$-module and $N$ is an Artinian $R$-module. We define the concept of filter coregular sequence to determine the infimum of integers $i$ such that the generalized local homology $\textrm{H}^{\mathfrak{a}}_i(M, N)$ is not finitely generated as an $\widehat{R}^{\mathfrak{a}}$-module, where $\widehat{R}^{\mathfrak{a}}$ denotes the $\mathfrak{a}$-adic completion of $R$. In particular, if $R$ is a complete semi-local ring, then $\textrm{H}^{\mathfrak{a}}_i(M, N)$ is a finitely generated $\widehat{R}^{\mathfrak{a}}$-module for all non-negative integers $i$ if and only if $(0:_N\mathfrak{a}+\textrm{Ann}(M))$ has finite length.
\end{abstract}

\maketitle
\section{Introduction}

In this paper, we consider a commutative Noetherian ring $R$ with non-zero identity, and an ideal $\fa\subseteq R$, as well as two $R$-modules $M$ and $N$. We denote the $\fa$-adic completion of $N$ by $\Lambda_\fa(N)$, and note that the $\fa$-adic completion functor $\Lambda_\fa(\cdot)$ is an additive covariant functor on the category of $R$-modules. We use $L_i^\fa(\cdot)$ to denote the $i$-th left derived functor of $\Lambda_\fa(\cdot)$. However, since the tensor functor is not left exact and the inverse limit is not right exact on the category of $R$-modules, computing the left-derived functors of $\Lambda_\fa(\cdot)$ is generally difficult. Moreover, it is important to note that $L_0^{\fa}(\cdot)\ncong\Lambda_\fa(\cdot)$.

Matlis studied $L^\fa_i(\cdot)$ in the case where $\fa$ is generated by a regular sequence and $R$ is a local ring in \cite{m1, m2}, and proved some duality between this functor and the local cohomology functor. Recently, Divaani-Aazar et al. in \cite{dft} studied the containment of $L^\fa_i(\cdot)$ in a Serre class of $R$-modules up to a given upper bound $s\geq 0$.

Cuong and Nam in \cite{cn1} defined the $i$-th {\it local homology} $\h {\fa}iN$ of $N$ with respect to $\fa$ as follows:
$$\h {\fa}iN:=\underset{n\in\N}{\varprojlim}\tor Ri{R/\fa^n}N.$$
They also showed that $L^{\fa}_i(N)\cong\h {\fa}iN$ when $N$ is Artinian. Similarly, the $i$-th {\it generalized local homology} $\gh {\fa}iMN$ of $M$ and $N$ with respect to $\fa$ is defined by
$$\gh {\fa}iMN:=\underset{n\in\N}{\varprojlim}\tor Ri{M/\fa^nM}N;$$
see \cite{n1} for basic properties and more details.

Matlis in \cite{m0} introduced the concept of cosequence (or coregular sequence) as a dual of the concept of regular sequence (see \cite{o} and \cite{tz} for more details and basic properties). If $N$ is Artinian and $(0:_N\fa)\neq 0$, then all maximal coregular $N$-sequences in $\fa$ have the same length, denoted by $\width {\fa}N$, where  $(0:_N\fa)$ denotes the set of all elements $x\in N$ such that $rx=0$ for all $r\in\fa$. Moreover, $\width {\fa}N=\inf\{i\in\mathbb Z: \h {\fa}iN\neq 0\}$ (see \cite[Theorem  4.11]{cn2}).

The  filter regular sequences can be used to study the Artinianess of local cohomology modules of finitely generated $R$-modules (see \cite[Sec. 3]{f}).
In this paper as a dual of the concept of filter regular sequence, we introduce the concept of filter coregular sequence   to study the finiteness of local homology modules of Artinian $R$-modules.

Let $\cosupp(N)$ denote the set of all prime ideals of $R$ containing $\ann(N)$.
A sequence $x_1,\dots,x_n$ of elements of $\fa$ is called a filter coregular $N$-sequence (of length $n$) in $\fa$ if $$\cosupp\big((0:_N(x_1,\dots,x_{i-1})R)/x_i(0:_N(x_1,\dots,x_{i-1})R)\big)\subseteq \Max (R)$$
for all $1\leq i\leq n$, where $\Max(R)$ denotes the set of all maximal ideals of $R$.

Assuming that $M$ is finitely generated and $N$ is Artinian, we prove that if there exists a filter coregular $N$-sequence in $\fa$ of infinite length, then every filter coregular $N$-sequence in $\fa$ can be extended to a filter coregular $N$-sequence in $\fa$ of infinite length, and in this case we set $\fwidth{\fa}N=\infty$. Now suppose that all filter coregular $N$-sequences in $\fa$ have finite length. Then all maximal filter coregular $N$-sequences in $\fa$ are of the same length, denoted by $\fwidth {\fa}N$. We prove (see Theorem \ref{thm} and Remark \ref{rem}) that:
\begin{align*}
  &\fwidth {\ann(M)}N\\
  &=\inf\{i\in\N_0: \tor RiMN\textrm{ has infinite length as an } R\textrm{-module}\}
  \end{align*}
and
 \begin{align*}
 &\fwidth {\fa+\ann(M)}N\\
  &=\inf\{i\in\N_0: \gh {\fa}iMN\textrm{ is not a finitely generated } \widehat{R}^{\fa}\textrm{-module}\}.
  \end{align*}
In particular,
 \begin{align*}
  &\fwidth {\fa}N\\
  &=\inf\{i\in\N_0: \h {\fa}iN\textrm{ is not a finitely generated } \widehat{R}^{\fa}\textrm{-module}\}.
  \end{align*}
We also show in Corollary \ref{cor} that if $\gh {\fa}iMN$ is a finitely generated $\widehat{R}^\fa$-module for all $i\in\N_0$, then $(0:_N\fa+\ann(M))$ has finite length. The converse statement is true when $R$ is a semi-local ring that is complete with respect to its Jacobson radical.

\section{Main results}

We shall use the following notations and terminologies. Let $\fa$ be an ideal of $R$ and $N$ be an $R$-module.  The radical of $\fa$ will be denoted by $\sqrt{\fa}$; also, $\ann(N)$ will denote the ideal $\{r\in R: rx=0 \textrm{ for all } x\in N\}$ of $R$; and $(0:_N\fa)$ will denote  the submodule $\{x\in N: rx=0 \textrm{ for all } r\in \fa\}$ of $N$. We denote by $\V{\fa}$ the set of all prime ideals of $R$ containing $\fa$; and  we use $\cosupp(M)$ to denote $\V{\ann(M)}$. The symbol $\mathbb N$ (respectively $\N_0$) will denote the set of positive (respectively non-negative) integers. We refer the reader for any unexplained terminology or notation to \cite{bs, mat, r}.

\begin{definition} Let $N$ be an $R$-module. We say a prime ideal $\fp$ of $R$ is   an {\it attached prime} of  $N$, if there exists a submodule $M$ of $N$ such that $\fp=\ann(N/M)$. We denote by $\att(N)$ the set of all attached primes of $N$.
\end{definition}
For an $R$-module $N$, it is clear that $\att(N)\subseteq\cosupp(N)$ (we refer the reader to \cite{o} for basic properties and more details of these notations).
When $N$ has a {\it secondary representation} in the sense of \cite{m}, our definition of $\att (N)$ coincides with that of Macdonald (see \cite[Exercise 7.2.5]{bs}). In particular, the set of attached primes of an Artinian module is a  finite set.

\begin{definition}
 Let $N$ be an $R$-module.  A sequence $x_1,\dots,x_n$ of elements of $R$ is called a filter coregular $N$-sequence (of length $n$) whenever $$\cosupp\big((0:_N(x_1,\dots,x_{i-1})R)/x_i(0:_N(x_1,\dots,x_{i-1})R)\big)\subseteq \Max (R)$$
for all $1\leq i\leq n$, where $\Max(R)$ denotes the set of all maximal ideals of $R$. If, in addition, $x_1,\dots,x_n$ belong to an ideal $\fa$, then we say that $x_1,\dots,x_n$ is a filter coregular $N$-sequence in $\fa$.
\end{definition}

\begin{lemma}\label{att}
Let $N$ be an Artinian $R$-module. The following conditions are equivalent:
\begin{enumerate}[\rm(i)]
  \item  $N$ has finite length;
  \item $\cosupp(N)\subseteq\Max(R)$; and
  \item   $\att(N)\subseteq \Max(R)$.
\end{enumerate}
\end{lemma}
\begin{proof}
 Assume that $N$ has finite length. Since $N$ is  finitely generated, we have $\cosupp(N)=\supp(N)$. Also,  the Artinianness of $N$ implies  that $\supp(N)\subseteq \Max(R)$, and so  $\cosupp(N)\subseteq\Max(R)$. This proves the implication (i)$\Rightarrow$(ii). The implication (ii)$\Rightarrow$(iii) is clear.
Finally, to prove the implication (iii)$\Rightarrow$(i), suppose that $\att(N)\subseteq \Max(R)$.  Then, by \cite[Proposition 7.2.11]{bs}, we have
$\sqrt{\ann(N)}=\bigcap_{\fq\in\att(N)}\fq$,
and so  $\big(\bigcap_{\fq\in\att(N)}\fq\big)^nN=0$ for some positive integer $n$.
It follows that $N$ has finite length because  $\att(N)$ consists of  finitely many maximal ideals.
\end{proof}

\begin{proposition}
Let  $x_1,\dots,x_n$  be elements of $R$,  and let $N$  be an Artinian $R$-module. The following conditions are equivalent:
\begin{enumerate}[\rm(i)]
\item $x_1,\dots,x_n$ is a filter coregular $N$-sequence;
\item $(0:_N(x_1,\dots,x_{i-1})R)/x_i(0:_N(x_1,\dots,x_{i-1})R)$ has finite length for all $1\leq i\leq n$;
\item $\att\left((0:_N(x_1,\dots,x_{i-1})R)/x_i(0:_N(x_1,\dots,x_{i-1})R)\right)\subseteq \Max (R)$ for all $1\leq i\leq n$; and
\item $x_i\notin\bigcup_{\fp\in\att\left(0:_N(x_1,\dots,x_{i-1})R\right)\setminus\Max (R)}\fp$ for all $1\leq i\leq n$.
\end{enumerate}
\end{proposition}
\begin{proof}
The statements (i)--(iii) are equivalent by Lemma \ref{att}. For each $1\leq i\leq n$, we set $N_{i-1}:=(0:_N(x_1,\dots,x_{i-1})R)$.   Then, in view of \cite[Proposition 2.13]{o}, we have  $x_i\notin\bigcup_{\fp\in\att(N_{i-1})\setminus\Max (R)}\fp$ if and only if $\att(N_{i-1}/x_i N_{i-1})=\V{x_iR}\cap\att(N_{i-1})\subseteq\Max (R)$.  Therefore (iii) and (iv) are  also equivalent.
\end{proof}

\begin{proposition}\label{prop} Let $M$ and $N$ be  $R$-modules, and let $x_1,\dots,x_n$ be elements of $R$. For each $i\in\N_0$, there are the following inclusions:
\begin{align}\label{1}
&\cosupp\left(\tor RiM{(0:_N(x_1,\dots,x_n)R)}\right)\subseteq\\
&\nonumber\left(\bigcup_{j=i}^{i+n} \cosupp\left(\tor RjMN\right)\right)\cup&\\
&\nonumber\left(\bigcup_{k=1}^n\bigcup_{j=i+2}^{i+2+n-k} \cosupp\left(\tor RjM{\frac{(0:_N(x_1,\dots,x_{k-1})R)}{x_k(0:_N(x_1,\dots,x_{k-1})R)}}\right)\right);
 \end{align}
and if, in addition,  $x_1,\dots,x_n$ belong to $\ann(M)$, then
\begin{align}\label{2}
&\cosupp\left(\tor RiMN\right)\subseteq \\
&\nonumber\cosupp\left(\tor R{i-n}M{(0:_N(x_1,\dots,x_n)R)}\right)\cup\\
&\nonumber\left(\bigcup_{k=1}^{n}\ \bigcup_{j=i+1-k}^{i+2-k}\cosupp\left(\tor RjM{\frac{(0:_N(x_1,\dots,x_{k-1})R)}{x_k(0:_N(x_1,\dots,x_{k-1})R)}}\right)\right).
\end{align}
\end{proposition}
\begin{proof}
We prove the claimed inclusions by induction on $n$. The following commutative diagram with exact rows
  \begin{align*} \xymatrix@C=.4cm{
 0\ar[r]&(0:_Nx_1R)\ar[r]&N\ar[r]^{x_1}\ar[d]^{x_1}\ar[dr]^{x_1}&x_1N\ar[r]&0\\
 &0\ar[r]&x_1N\ar[r]^{\subseteq}&N\ar[r]&N/x_1N\ar[r]&0
} \end{align*}
 induces the  commutative diagram
  \begin{align}\label{3}\xymatrix@C=0.4cm{
 \cdots\ar[r] &T_i(0:_Nx_1R)\ar[r]&T_i(N)\ar[r]^{x_1^{(i)}}\ar[d]^{x_1^{(i)}}\ar[rd]^{x_1}&T_i(x_1N)\ar[r]&T_{i-1}(0:_Nx_1R)\ar[r]&\cdots\\
 \cdots\ar[r]&T_{i+1}(N/x_1N)\ar[r]&T_i(x_1N)\ar[r]^{f_i}&T_i(N)\ar[r]&T_i(N/x_1N)\ar[r]&\cdots
} \end{align}
with exact rows, where $T_i(\cdot):=\tor RiM{\cdot}$ and $x_1^{(i)}:=\tor Ri{{\rm id}_M}{x_1}$. Therefore \cite[Proposition 2.9(4)]{o} implies that
 \begin{align}\label{4}
&\cosupp\left(T_i(0:_Nx_1R)\right)\subseteq
\cosupp\left(T_i(N)\right)\cup\cosupp\left(T_{i+1}(x_1N)\right)\\
 &\nonumber\subseteq\cosupp\left(T_i(N)\right)\cup\cosupp\left(T_{i+1}(N)\right)\cup\cosupp\left(T_{i+2}(N/x_1N)\right)
 \end{align}
for all $i\in\N_0$ (we note that if $L\rightarrow M\rightarrow N$ is an exact sequence of $R$-modules, then we can deduce from \cite[Proposition 2.9(4)]{o} that $\cosupp(M)\subseteq\cosupp(L)\cup\cosupp(N)$). This proves   the inclusion (\ref{1}) in the case when $n=1$. Now assume, inductively, that $n>1$ and the inclusion (\ref{1}) holds for  smaller values of $n$. If we replace $N$ by $(0:_Nx_1R)$, then, by the inductive hypothesis for elements $x_2,\dots,x_n$, we have
 \begin{align}\label{5}
&\cosupp\left(T_i{(0:_N(x_1,\dots,x_n)R)}\right)\subseteq\left(\bigcup_{j=i}^{i+n-1} \cosupp\left(T_j{(0:_Nx_1R)}\right)\right)\\
&\nonumber\cup \left( \bigcup_{k=2}^{n}\bigcup_{j=i+2}^{i+2+n-k} \cosupp\left(T_j\left(\frac{(0:_N(x_1,\dots,x_{k-1})R)}{x_k(0:_N(x_1,\dots,x_{k-1})R)}\right)\right)\right)
 \end{align}
 (note that if we  set $y_1:=x_2,\dots,y_{n-1}:=x_n$ and $l:=k-1$, then $1\leq l\leq n-1$ and $i+2\leq j\leq i+2+n-1-l$  yield  $2\leq k\leq n$ and $i+2\leq j\leq i+2+n-k$).
Now combining the inclusion (\ref{4})  with the inclusion (\ref{5}) yields the inclusion  (\ref{1}) and  the inductive step is complete.

  Now assume that $x_jM=0$ for all $1\leq j\leq n$ and we prove, by induction on $n$, that the inclusion (\ref{2}) holds.   Since the functor $T_i(\cdot)$ is $R$-linear, the endomorphism of $T_i(N)$ given by multiplication by $x_j$  is the zero map for all $i\in\N_0$ and all $1\leq j\leq n$. The triangle in  the  diagram $(\ref{3})$ commutes, and so  $\im x_1^{(i)}\subseteq\Ker f_i$ for all $i\in\N_0$. Therefore 
\begin{equation}\label{6}\cosupp\left(\im x^{(i)}_1\right)\subseteq\cosupp\left(\Ker f_i\right)\subseteq\cosupp\left(T_{i+1}(N/x_1N)\right).\end{equation}
 Also, the exactness of rows in the diagram (\ref{3}) implies that
\begin{equation}\label{7}\cosupp\left(T_i(x_1N)\right)\subseteq\cosupp\left(\im x^{(i)}_1\right)\cup\cosupp\left(T_{i-1}(0:_Nx_1)\right)\end{equation}
and
\begin{equation}\label{8}\cosupp\left(T_i(N)\right)\subseteq\cosupp\left(T_i(x_1N)\right)\cup\cosupp\left(T_i(N/x_1N)\right).\end{equation}
The inclusions (\ref{6}), (\ref{7}) and (\ref{8}) yield
\begin{align}\label{9}
&\cosupp\left(T_i(N)\right)\subseteq\cosupp\left(T_{i-1}(0:_Nx_1R)\right)\cup\\
&\nonumber\cosupp\left(T_i\left({N}/{x_1N}\right)\right)\cup\cosupp\left(T_{i+1}\left({N}/{x_1N}\right)\right).
\end{align}
Hence,  the inclusion (\ref{2}) is true in the case when $n=1$. Next suppose, inductively,  that $n>1$ and that the inclusion (\ref{2}) has been proved for smaller values of $n$.  If we use $(0:_Nx_1R)$ and $i-1$ instead of $N$ and $i$  respectively, then the inductive hypothesis for elements $x_2,\dots,x_n$ yields
\begin{align}\label{10}
&\cosupp\left(T_{i-1}(0:_Nx_1R)\right)\subseteq \cosupp \left(T_{i-n}{(0:_N(x_1,\dots,x_n)R)}\right)\\
&\nonumber\cup\left(\bigcup_{k=2}^{n}\ \bigcup_{j=i+1-k}^{i+2-k}\cosupp
\left(T_j\left({\frac{(0:_N(x_1,\dots,x_{k-1})R)}{x_k(0:_N(x_1,\dots,x_{k-1})R)}}\right)\right)\right).
\end{align}
By combining the inclusions  (\ref{9}) and (\ref{10}), we obtain the inclusion (\ref{2}). This completes  the inductive step.
\end{proof}

\begin{corollary}\label{corcosupp} Let $M$ and $N$ be $R$-modules, and let $x_1,\dots,x_n$ be a filter coregular $N$-sequence in $\ann(M)$. Then
 \begin{align}\label{11}
 \cosupp\left(\tor RiMN\right)\subseteq\Max(R)\end{align}
 for all $i<n$, and
 \begin{align}\label{12}
 &\cosupp\left(\tor RnMN\right)\cup\Max(R)=\\
 &\nonumber\cosupp\left(M\otimes_R(0:_N(x_1,\dots,x_n)R)\right)\cup\Max(R).\end{align}
\end{corollary}
\begin{proof}
For each $1\leq k\leq n$, since $\cosupp\left(\frac{(0:_N(x_1,\dots,x_{k-1})R)}{x_k(0:_N(x_1,\dots,x_{k-1})R)}\right)\subseteq\Max(R)$, we have
\begin{equation}\label{13}
\cosupp\left(\tor RiM{\frac{(0:_N(x_1,\dots,x_{k-1})R)}{x_k(0:_N(x_1,\dots,x_{k-1})R)}}\right)\subseteq\Max(R)
\end{equation}
 for all $i\in\N_0$.
Hence the inclusion (\ref{11}) is an immediate consequence of the inclusion (\ref{2}).
Now we prove  the equation (\ref{12}). If we set $i=0$ in the inclusion (\ref{1}), then it follows from the inclusions (\ref{11})  and (\ref{13}) that
\begin{align}\label{14}
&\cosupp\left(M\otimes_R(0:_N(x_1,\dots,x_n)R)\right)\subseteq \\
&\nonumber\cosupp\left(\tor RnMN\right)\cup\Max(R).
\end{align}
Conversely, if we set $i=n$ in the inclusion (\ref{2}), then the inclusion (\ref{13}) implies that
\begin{align}\label{15}
&\cosupp\left(\tor RnMN\right)\subseteq\\
&\nonumber \cosupp\left(M\otimes_R(0:_N(x_1,\dots,x_n)R)\right)\cup\Max(R).
\end{align}
Now the equation (\ref{12}) follows from the inclusions (\ref{14}) and (\ref{15}).
\end{proof}
\begin{lemma}\label{lemnew}
Let  $M, N$ and $L$ be  $R$-modules such that   $M$ and $L$ are finitely generated, and let $n\in\N$. If $\cosupp \left(\tor RiMN\right)\subseteq\Max(R)$ for all $i<n$ and $\supp(L)\subseteq \supp(M)$, then
$\cosupp \left(\tor RiLN\right)\subseteq\Max(R)$ for all $i<n$. In particular, $\cosupp \left(\tor RiLN\right)\subseteq\Max(R)$ for all $i<n$ if and only if $\cosupp \left(\tor RiMN\right)\subseteq\Max(R)$ for all $i<n$
whenever $\supp(L)=\supp(M)$.
\end{lemma}
\begin{proof}
Assume that $\cosupp \left(\tor RiMN\right)\subseteq\Max(R)$ for all $i<n$ and we prove by induction on $n$ that  for every finitely generated $R$-module $L$ with $\supp(L)\subseteq \supp(M)$, $\cosupp \left(\tor RiLN\right)\subseteq\Max(R)$ for all $i<n$. Since $\supp(L)\subseteq \supp(M)$, by Gruson's theorem \cite[Theorem 4.1]{v} there exists a chain
$$0=L_0\subseteq L_1\subseteq\dots\subseteq L_m=L$$
of submodules of $L$ such that, for each $1\leq j\leq m$,  $L_j/L_{j-1}$ is a homomorphic image of a direct sum of finitely many copies of $M$.
For each $1\leq j\leq m$, the exact sequence
\begin{equation*}
0\rightarrow L_{j-1}\rightarrow L_j\rightarrow L_j/L_{j-1}\rightarrow 0
\end{equation*}
induces the following exact sequence
\begin{align*}
 \dots\rightarrow \tor Ri{L_{j-1}}N\rightarrow \tor Ri{L_j}N\rightarrow \tor Ri{L_j/L_{j-1}}N\rightarrow\cdots.
 \end{align*}
 Hence
\begin{align*}
&\cosupp(\tor Ri{L_j}N)\\
&\subseteq\cosupp(\tor Ri{L_{j-1}}N)\cup\cosupp(\tor Ri{L_j/L_{j-1}}N)\\
\end{align*}
for all $1\leq j\leq m$ and all $i$. It follows that
\begin{align*}
&\cosupp(\tor Ri{L}N)\\
&=\cosupp(\tor Ri{L_m}N)\\
&\subseteq\cosupp(\tor Ri{L_{m-1}}N)\cup\cosupp(\tor Ri{L_m/L_{m-1}}N)\\
&\ \vdots\\
&\subseteq\cosupp(\tor Ri{L_{0}}N)\cup\left(\bigcup_{j=1}^m\cosupp(\tor Ri{L_j/L_{j-1}}N\right)\\
&=\bigcup_{j=1}^m\cosupp(\tor Ri{L_j/L_{j-1}}N
\end{align*}
for all $i<n$. Thus to prove the assertion it is sufficient for us to prove that $\cosupp(\tor Ri{L_j/L_{j-1}}N\subseteq\Max(R)$
for all $1\leq j\leq m$ and all $i<n$. Hence we can assume that $m=1$ and there exists an exact sequence
\begin{equation*}
0\rightarrow K\rightarrow M^t\rightarrow L\rightarrow 0
\end{equation*}
 for some $t\in\N$ and some finitely generated $R$-module $K$. This exact sequence induces the following long exact sequence
 \begin{align}\label{new2}
 \dots\rightarrow \tor RiMN^{t}\rightarrow \tor RiLN\rightarrow \tor R{i-1}KN\rightarrow\cdots.
 \end{align}

 For $n=1$, it follows  from the exact sequence
$(M\otimes_RN)^{t}\rightarrow L\otimes_RN\rightarrow 0$ that
$$\cosupp(L\otimes_RN)\subseteq\cosupp(M\otimes_RN)\subseteq\Max(R).$$
Therefore the result holds for $n=1$. Now assume, inductively, that $n>1$ and the result has been proved for smaller values of $n$. It follows from  the exact sequence (\ref{new2}) that \begin{align}\label{new3}
&\cosupp(\tor RiLN)\\
&\nonumber\subseteq\cosupp(\tor Ri{M}N)\cup\cosupp(\tor R{i-1}KN)
\end{align}
 for all $i$. Since $\supp(K)\subseteq\supp(M)$, the induction hypothesis implies that $$\cosupp(\tor RiKN)\subseteq\Max(R)$$ for all $i<n-1$.
Thus, by the hypothesis and the inclusion (\ref{new3}), we have $$\cosupp(\tor RiLN)\subseteq\Max(R)$$ for all $i<n$.  This completes the inductive step.
\end{proof}

 Now, we are ready to state and prove the main result of this paper. Let $\fa$ be an ideal of $R$ and let $N$ be an Artinian $R$-module. Among the other things, the following theorem shows that  the infimum of integers $i$ with the property  that the  local homology module $\h {\fa}iN$  is not finitely generated as an $\widehat{R}^\fa$-module and the common length of all maximal  filter coregular $N$-sequences in $\fa$ are same.
\begin{theorem}\label{thm} Let $\fa$ be an ideal of $R$, and  let  $M$ and $N$  be $R$-modules such that $M$ is finitely generated and $N$ is Artinian. For each $n\in\N$,  the following conditions are equivalent:
\begin{enumerate}[\rm(i)]
\item there is a  filter coregular $N$-sequence in $\fa$ of length $n$;
\item any filter coregular $N$-sequence in $\fa$ of length less than $n$ can be extended to a filter coregular $N$-sequence in $\fa$ of length $n$;
 \item $\cosupp\left(\tor Ri{R/\fa}N\right)\subseteq\Max (R)$ (or equivalently $\tor Ri{R/\fa}N$ has finite length)  for all $i<n$;
\item if $\supp(M)=\V{\fa}$, then  $\cosupp\left(\tor RiMN\right)\subseteq\Max (R)$ (or equivalently $\tor Ri{M}N$ has finite length)  for all $i<n$; and
\item if $\ann(M)\subseteq\fa$,  then $\gh {\fa}iMN$  is a finitely generated $\widehat{R}^\fa$-module for all $i<n$.
\end{enumerate}
\end{theorem}
\begin{proof}
The statements (iii) and (iv) are equivalent by Lemma \ref{lemnew}.  The implication (ii)$\Rightarrow$(i) is clear.  Also, (i)$\Rightarrow$(iii) is an immediate consequence of the inclusion (\ref{11}) in  Corollary \ref{corcosupp}.

(iii)$\Rightarrow$(ii). Assume that  $\cosupp\left(\tor Ri{R/\fa}N\right)\subseteq\Max (R)$ for all $i<n$, and suppose, for the sake of contradiction, that  $x_1,\dots,x_m$ is a maximal filter coregular $N$-sequence  in $\fa$ of length $0\leq m<n$. The  maximality of  $x_1,\dots,x_m$ yields
 $$\fa\subseteq\bigcup_{\fp\in\att\left(0:_N(x_1,\dots,x_m)R\right)\setminus\Max (R)}\fp.$$
Since $\att\left(0:_N(x_1,\dots,x_m)R\right)$ is a finite set, it follows from the Prime Avoidance Theorem that  $\fa\subseteq\fp$ for some $\fp\in\att\left(0:_N(x_1,\dots,x_m)R\right)\setminus\Max (R)$. Hence,  by the equation (\ref{12}) in Corollary \ref{corcosupp} and the hypothesis, we have
\begin{align*}
\fp\in&\V{\fa}\cap\att(0:_N(x_1,\dots,x_m)R)\\
&=\att(R/\fa\otimes_R(0:_N(x_1,\dots,x_m)R))\\
&\subseteq\cosupp(R/\fa\otimes_R(0:_N(x_1,\dots,x_m)R))\\
&\subseteq\cosupp(\tor Rm{R/\fa}N)\cup\Max(R)\\
&\subseteq\Max(R),
\end{align*}
which is a contradiction.  Hence the statements (i)--(iv)  are equivalent.

(i)$\Leftrightarrow$(v). We prove, by induction on $n$,  that (i) and (v) are equivalent. Assume that $M$ is a finitely generated $R$-module such that $\ann(M)\subseteq\fa$.
 We first assume that $n=1$.  Since $M\otimes_RN$ is Artinian, we have
\begin{align*}
\gh {\fa}0MN\cong \Lambda_\fa(M\otimes_RN)\cong (M\otimes_R N)/\fa^s(M\otimes_R N)
\cong \tor R0{M/\fa^sM}N
\end{align*}
for all sufficiently large integers $s$. Also, since $\supp(M/\fa^sM)=\V\fa$,  the equivalence of (i) and (iv)
implies that $\gh {\fa}0MN$ is a finitely generated $R$-module or equivalently it is a finitely generated $\widehat{R}^{\fa}$-module if and only if $\fa$ contains a filter coregular element on $N$ (note that since $\gh {\fa}0MN$ is $\fa$-torsion, its submodules as an $R$-module and as an $\widehat{R}^{\fa}$-module are same; see  \cite[Lemma 1.3]{kls}). Thus the result holds in the case $n=1$.

Now assume, inductively, that $n>1$ and  the result has been proved for smaller values of $n$. Since $N$ is Artinian, there exists $t\in\N$ such that $\fa^sN=\fa^tN$ for all $s\geq t$ and so $\Lambda_\fa(N)\cong N/\fa^tN$. Assume that either (i) or (v) holds. Since $n>1$ and $\gh {\fa}0MN\cong\tor R0{M/\fa^sM}N$ for sufficiently large integers $s$,  if (v) holds, then $\tor R0{M/\fa^sM}N$ has finite length by the hypothesis of (v). Since $\supp(R/\fa^t)=\supp(M/\fa^sM)$, by Lemma \ref{lemnew},  $\Lambda_\fa(N)\cong\tor R0{R/\fa^t}N$ has finite length in this case. Also, if (i) holds, then, by the equivalence of (i) and (iv), $\Lambda_\fa(N)\cong\tor R0{R/\fa^t}N$ has finite length. Therefore in the either cases $\Lambda_\fa(N)$ has finite length.
  Now, the exact sequence
$$0\rightarrow \fa^tN\rightarrow N\rightarrow \Lambda_\fa(N)\rightarrow 0$$
of Artinian $R$-modules induces the following long exact sequences 
\begin{align}\label{18}
\cdots&\rightarrow\gh {\fa}{i+1}M{\Lambda_\fa(N)}\rightarrow\gh {\fa}iM{\fa^tN}\rightarrow\gh {\fa}iMN\\
&\nonumber\rightarrow\gh {\fa}iM{\Lambda_\fa(N)}\rightarrow\cdots
\end{align}
(see \cite[Proposition 2.4]{n1}), and
\begin{align}\label{19}
\cdots&\rightarrow\tor R{i+1}{R/\fa}{\Lambda_\fa(N)}\rightarrow\tor R{i}{R/\fa}{\fa^tN}\rightarrow\tor R{i}{R/\fa}N\\
&\nonumber\rightarrow\tor Ri{R/\fa}{\Lambda_\fa(N)}\rightarrow\cdots.
\end{align}
Since   $\Lambda_\fa(N)$ is Artinian, by  \cite[Theorems 2.3(i) and 3.2]{n2}, we have
$$\gh {\fa}iM{\Lambda_\fa(N)}\cong H_i\left(\Lambda_\fa(M\otimes_RF_{\bullet}) \right),$$
where $F_\bullet$ is a free resolution of $\Lambda_\fa(N)$. Now $\Lambda_\fa(N)$ is finitely generated and so we can assume that every component of $F_\bullet$ is finitely generated. On the other hand, $\Lambda_\fa(\cdot)$ is an additive exact functor on the category of finitely generated $R$-modules, and hence it commutes with the homological functor in this category. Therefore
$$ H_i\left(\Lambda_\fa(M\otimes_RF_{\bullet}) \right)\cong\Lambda_\fa\left(H_i\left(M\otimes_RF_{\bullet}\right)\right)\cong \Lambda_\fa\left(\tor RiM{\Lambda_\fa(N)}\right).$$
Since $\tor RiM{\Lambda_\fa(N)}$ is Artinian, we obtain
$$\Lambda_\fa\left(\tor RiM{\Lambda_\fa(N)}\right)\cong\tor RiM{\Lambda_\fa(N)}/\fa^r\tor RiM{\Lambda_\fa(N)}$$ for all sufficiently large integers $r$.
Since $\fa^r\tor RiM{\Lambda_\fa(N)}=0$ for $r\geq t$, the above isomorphisms yield $$\gh {\fa}iM{\Lambda_\fa(N)}\cong\tor RiM{\Lambda_\fa(N)}\cong\tor RiM{N/\fa^tN}.$$
Hence $\gh {\fa}iM{\Lambda_\fa(N)}$
 is a finitely generated $R$-module  for all $i\in\N_0$. Also,  the above isomorphism shows that $\gh {\fa}iM{\Lambda_\fa(N)}$ is $\fa$-torsion, and so
 $\gh {\fa}iM{\Lambda_\fa(N)}$ is a finitely generated  $\widehat {R}^\fa$-module  for all $i\in\N_0$ by \cite[Lemma 1.3]{kls}. Now, for each $i\in\N_0$, it follows from the long exact sequence (\ref{18}) that $\gh {\fa}iMN$ is a finitely generated $\widehat {R}^\fa$-module if and only if $\gh {\fa}iM{\fa^tN}$ is a finitely generated $\widehat {R}^\fa$-module. Also, for each $i\in\N_0$,
it follows from the long exact sequence (\ref{19}) that $\tor Ri{R/\fa}N$ has finite length if and only if $\tor Ri{R/\fa}{\fa^tN}$ has finite length because  $\tor RiM{\Lambda_\fa(N)}$ has finite length for all $i$. Thus to prove the equivalence of (i) and (v), in view of the equivalence of (i) and (iii),  we can replace $N$ by
   $\fa^tN$ and assume, in addition, that $\fa N=N$. Therefore $\V{\fa}\cap\att(N)=\emptyset$, and so
   $\fa\nsubseteq\bigcup_{\fp\in\att(N)}\fp$. Let $x_1\in\fa\setminus\bigcup_{\fp\in\att(N)}\fp$.
    Then $\V{x_1R}\cap\att(N)=\emptyset$, and so $N=x_1N$.
The exact sequence
$$0\rightarrow (0:_Nx_1R)\rightarrow N\stackrel{x_1}{\longrightarrow}N\rightarrow 0$$
induces the long exact sequence
\begin{align}\label{20}
\cdots&\rightarrow\gh {\fa}{i+1}MN\stackrel{x_1}{\longrightarrow}\gh {\fa}{i+1}MN\rightarrow\gh {\fa}iM{(0:_Nx_1R)}\\
&\nonumber\rightarrow\gh {\fa}{i}MN\stackrel{x_1}{\longrightarrow}\cdots.
\end{align}
We first assume that  (i) holds. By the equivalence of (i) and (ii), we can extend $x_1$ to a filter coregular $N$-sequence  of length $n$, say $x_1,x_2,\dots,x_n$. Hence $x_2,\dots,x_n$ is a filter coregular $(0:_Nx_1R)$-sequence in $\fa$, and so, by the inductive hypothesis, $\gh {\fa}iM{(0:_Nx_1R)}$ is a finitely generated $\widehat{R}^{\fa}$-module for all $i<n-1$. It follows from the  long exact sequence (\ref{20}) that $\gh {\fa}iMN/x_1\gh {\fa}iMN$ and consequently its homomorphic image   $\gh {\fa}iMN/(\fa\widehat{R}^\fa)\gh {\fa}iMN$ are finitely generated $\widehat{R}^{\fa}$-modules for all $i<n$. Also, by \cite[Proposition 2.3(i)]{n1}, we have
$$\bigcap_{t\in\N}(\fa\widehat{R}^\fa)^t\gh {\fa}iMN=\bigcap_{t\in\N}({\fa}^t\widehat{R}^\fa)\gh {\fa}iMN=\bigcap_{t\in\N}{\fa}^t\gh {\fa}iMN=0.$$
   Hence, by  \cite[Theorem 8.4]{mat}, $\gh {\fa}iMN$ is a finitely generated $\widehat{R}^{\fa}$-module for all $i<n$. Conversely, assume that $\gh {\fa}iMN$ is a finitely generated $\widehat{R}^{\fa}$-module for all $i<n$. It follows from the  long exact sequence (\ref{20}) that $\gh {\fa}{i}M{(0:_Nx_1R)}$  is a finitely generated $\widehat{R}^{\fa}$-module for all $i<n-1$, and so, by the inductive hypothesis, there is a   filter coregular $(0:_Nx_1R)$-sequence in $\fa$ of length $n-1$,  say $x_2,\dots,x_n$. Therefore $x_1,x_2,\dots,x_n$ is a filter coregular $N$-sequence. This completes the inductive step.
\end{proof}

\begin{remark}\label{rem}  Let $\fa$ be an ideal of $R$, and let $N$ be an Artinian $R$-module. When there exists a filter coregular $N$-sequence in $\fa$ of infinite length, then, by the equivalence of  (i) and (ii) in Theorem \ref{thm},  any filter coregular $N$-sequence in $\fa$ can be extended to a filter coregular $N$-sequence in $\fa$ of arbitrary length, and in this case we set $\fwidth{\fa}N=\infty$. Now assume that all filter coregular $N$-sequences in $\fa$ have finite length.  Again, by the equivalence of (i) and (ii) in Theorem \ref{thm},
 we can extend any filter coregular $N$-sequence in $\fa$ to a maximal one, and all  maximal filter coregular $N$-sequences in $\fa$ are of the  same  length which we denote this common length by $\fwidth {\fa}N$.
     Moreover, if  $M$ is a finitely generated $R$-module such that $\supp(M)=\V{\fa}$, then, by Theorem \ref{thm}, we have
  \begin{align}\label{21}
&\fwidth {\fa}N\\
 \nonumber &=\inf\{i\in\N_0: \cosupp\left(\tor RiMN\right)\nsubseteq \Max (R)\}\\
 \nonumber &=\inf\{i\in\N_0: \tor RiMN\textrm{ has infinite length as an } R\textrm{-module}\}\\
   & \nonumber=\inf\{i\in\N_0: \h {\fa}iN\textrm{ is not a finitely generated } \widehat{R}^{\fa}\textrm{-module}\}
  \end{align}
(we note that $\gh {\fa}iRN=\h {\fa}iN$).
Also, for an arbitrary finitely generated $R$-module $L$, since $\gh {\fa+\ann(L)}iLN\cong\gh {\fa}iLN$, if we replace $\fa$  by $\fa+\ann(L)$ in Theorem \ref{thm},  then  the equivalence of (ii) and (v) in Theorem \ref{thm} yields
\begin{align}\label{22}
 &\fwidth {\fa+\ann(L)}N\\
   & \nonumber=\inf\{i\in\N_0: \gh {\fa}iLN\textrm{ is not a finitely generated } \widehat{R}^{\fa}\textrm{-module}\}.
  \end{align}
Finally,  since $\V{\fa}=\V{\sqrt\fa}$, it follows from the first equality in the equation  (\ref{21}) that $\fwidth{\fa}N=\fwidth{\sqrt\fa}N.$
\end{remark}
\begin{proposition}\label{infinity}
Let $\fa$ be an ideal of $R$, and  let $N$ be an Artinian $R$-module. If $\fwidth {\fa}N=\infty$, then $(0:_N\fa)$ has finite length. The converse  statement  holds whenever $R$ is a semi-local ring which is complete with respect to its Jacobson radical.
\end{proposition}
\begin{proof}
 Assume that $\fwidth {\fa}N=\infty$, and $x_1, x_2, x_3,\dots$ is a filter coregular $N$-sequence of infinite length in $\fa$. There is the following descending chain of submodules of $N$
$$(0:_Nx_1R)\supseteq(0:_N(x_1, x_2)R)\supseteq(0:_N(x_1, x_2, x_3)R)\supseteq\cdots.$$
Hence $(0:_N(x_1,\dots,x_{n-1})R)=(0:_N(x_1,\dots,x_n)R)$ for some $n\in\N$, and so $x_n(0:_N(x_1,\dots,x_{n-1})R)=0$. Thus
 $(0:_N(x_1,\dots,x_{n-1})R)$ has finite length because   $(0:_N(x_1,\dots,x_{n-1})R)/x_n(0:_N(x_1,\dots,x_{n-1})R)$ has finite length by definition. Hence $(0:_N\fa)\subseteq(0:_N(x_1,\dots,x_{n-1})R)$ has finite length. To prove the converse statement,
assume that $R$ is a complete semi-local ring and that $(0:_N\fa)$ has finite length. Hence
$$\cosupp(0:_N\fa)=\supp(0:_N\fa)\subseteq\Max(R).$$
On the other hand, for each $i\in\N_0$,  $\fa+\ann(N)\subseteq\ann(\tor Ri{R/\fa}N)$. Therefore, in view of \cite[Proposition 2.12]{o}, we have
 \begin{align*}
& \cosupp(\tor Ri{R/\fa}N)\subseteq\V{\fa+\ann(N)}\\
 &=\V{\fa}\cap\cosupp(N) =\cosupp(0:_N\fa)\subseteq\Max (R)
 \end{align*}
 for all $i\in\N_0$. Hence Theorem \ref{thm} implies that  $\fwidth {\fa}N=\infty$.
\end{proof}
\begin{corollary}\label{cor}  Let $\fa$ be an ideal of $R$, and let  $M$ and $N$ be $R$-modules such that $M$ is finitely generated and $N$ is Artinian.
\begin{enumerate}[\rm(i)]
 \item If $\tor RiMN$ has finite length for all $i\in\N_0$, then $(0:_N\ann(M))$ has finite length.
\item If $\gh {\fa}iMN$ is a finitely generated $\widehat{R}^\fa$-module  for all $i\in\N_0$, then $(0:_N\fa+\ann(M))$ has finite length. In particular, $(0:_N\fa)$ has finite length whenever $\h {\fa}iN$ is a finitely generated $\widehat{R}^\fa$-module  for all $i\in\N_0$.
 \end{enumerate}
  Moreover, the converse  statements  hold when $R$ is  a complete  semi-local ring.
\end{corollary}
\begin{proof}
It follows by the equations (\ref{21}), (\ref{22}) and Proposition \ref{infinity}.
\end{proof}


\vskip 0.4 true cm

\begin{center}{\textbf{Acknowledgments}}
\end{center}
The authors are deeply grateful to the referee for his/her useful comments and suggestions.  Also, we would like to express our gratitude to Prof. Kamran Divaani-Aazar for careful reading of the manuscript and for the helpful  comments and suggestions.
 \\ \\
\vskip 0.4 true cm
\bibliographystyle{amsplain}

\begin{thebibliography}{5}

\bibitem{bs}
M. P. Brodmann and R. Y. Sharp,
{\it Local Cohomology: An Algebraic Introduction with Geometric Applications},
Cambridge Studies in Advanced Mathematics {\bf 60}
 (Cambridge University Press, Cambridge, 1998).



\bibitem{cn1}
N. T. Cuong and T. T. Nam,
The I-adic completion and local homology for Artinian modules,
Math. Proc. Cambridge Philos. Soc. \textbf{131}(1) (2001), 61--72.

\bibitem{cn2}
N. T. Cuong and T. T. Nam,
A local homology theory for linearly compact modules,
J. Algebra \textbf{319}(11) (2008), 4712--4737.

\bibitem{dft}
K. Divaani-Aazar, H. Faridian and M. Tousi,
 Local homology, Koszul homology and Serre classes,
  Rocky Mountain J. Math.
  {\bf 48}(6) (2018) 1841--1869.

 \bibitem{f}
A. Fathi,
The first non-isomorphic local cohomology modules with respect to their ideals,
  J. Algebra Appl. {\bf17}(12) (2018) 1850230.

\bibitem{kls}
B. Kubik, M. J. Leamer and S. Sather-Wagstaff,
Homology of Artinian and Matlis reflexive modules, I,
J. Pure Appl. Algebra {\bf 215}(10) (2011), 2486--2503.


\bibitem{m}
I. G. Macdonald,
Secondary representation of modules over a commutative ring,
  Symposia Mathematica \textbf{11} (1973), 23--43.

\bibitem{m0}
E. Matlis, Modules with descending chain condition,
 Trans. Amer. Math. Soc. {\bf 97} (1960) 495--508.

\bibitem{m1}
E. Matlis,  The Koszul complex and duality,
Comm. Algebra {\bf 1} (1974), 87--144.

\bibitem{m2}
E. Matlis, The higher properties of $R$-sequences,
 J. Algebra {\bf 50}(1) (1978), 77--112.


\bibitem{mat}
H. Matsumura,
Commutative ring theory,
Cambridge Studies in Advanced Mathematics {\bf 8}
 (Cambridge University Press, Cambridge, 1986).

\bibitem{n2}
T. T. Nam, Left-derived functors of the generalized $I$-adic completion and generalized local homology,
Comm. Algebra \textbf{38}(2) (2010), 440--453.

\bibitem{n1}
T. T. Nam,
Generalized local homology for Artinian modules,
Algebra Colloq. \textbf{19}(1) (2012), 1205--1212.

\bibitem{o}
A. Ooishi,
Matlis duality and the width of a module,
Hiroshima Math. J.  \textbf{6}(3) (1976), 573--587.

\bibitem{r}
J. J. Rotman,
{\it An introduction to homological algebra}, Pure and Applied Mathematics {\bf 85} (Academic Press, Inc., New York--London, 1979).


\bibitem{tz}
Z. Tang and H. Zakeri,
Co-Cohen-Macaulay modules and modules of generalized fractions,
Comm. Algebra \textbf{22}(6) (1994) 2173--2204.

\bibitem{v}
W. V. Vasconcelos,
{\it Divisor Theory in Module Categories}, North-Holland Mathematics Studies {\bf 14} (North-Holland Publishing Co., Amsterdam, 1974).

\end{thebibliography}

\bigskip
\bigskip


{\footnotesize {\bf Ali Fathi}\; \\ {Department of
Mathematics}, {Zanjan Branch,
Islamic Azad University,} {Zanjan, Iran.}\\
{\tt Email:alif1387@gmail.com}\\

{\footnotesize {\bf Alireza Hajikarimi}\; \\ {Department of
Mathematics}, {Mobarakeh Branch,
Islamic Azad University,} {Isfahan, Iran.}\\
{\tt Email: a.hajikarimi@mau.ac.ir}\\

\end{document}